\pdfoutput=1
\documentclass[12pt,reqno]{amsart}

\usepackage[T1]{fontenc}
\usepackage[utf8]{inputenc}
\usepackage{lmodern}
\usepackage{microtype}
\usepackage{geometry}
\usepackage{amsmath,amssymb,amsthm,mathtools}
\usepackage{bm}
\usepackage{enumitem}
\usepackage{graphicx}
\usepackage{booktabs}
\usepackage{array}
\usepackage{xcolor}
\usepackage[colorlinks=true,linkcolor=blue,citecolor=teal,urlcolor=magenta]{hyperref}
\usepackage[nameinlink,capitalise,noabbrev]{cleveref}
\usepackage{tikz}
\usepackage{cite}

\numberwithin{equation}{section}

\theoremstyle{plain}
\newtheorem{theorem}{Theorem}[section]

\newtheorem{proposition}[theorem]{Proposition}

\theoremstyle{definition}

\theoremstyle{remark}
\newtheorem{remark}[theorem]{Remark}

\crefname{theorem}{Theorem}{Theorems}
\crefname{lemma}{Lemma}{Lemmas}
\crefname{proposition}{Proposition}{Propositions}
\crefname{corollary}{Corollary}{Corollaries}
\crefname{definition}{Definition}{Definitions}
\crefname{assumption}{Assumption}{Assumptions}
\crefname{remark}{Remark}{Remarks}
\crefname{example}{Example}{Examples}
\crefname{equation}{Equation}{Equations}
\crefname{section}{Section}{Sections}
\crefname{appendix}{Appendix}{Appendices}

\newif\ifnotes
\notestrue
\ifnotes
  \newcommand{\todo}[1]{\textcolor{red}{[TODO: #1]}}
\else
  \newcommand{\todo}[1]{}
\fi

\title[Eighth-Order Accurate Methods for Lane-Emden BVPs]{Eighth-Order Accurate Methods for Boundary Value Problems Arising from the Lane-Emden Equation}
\author[Dang Quang A]{Dang Quang A}
\address{Centre for Informatics and Computing, VAST,
	18 Hoang Quoc Viet, Cau Giay, Hanoi, Vietnam \\
	Thai Nguyen University of Information and Communication Technology,
	Thai Nguyen, Viet Nam }
\email{dangquanga@cic.vast.vn}

\author[Nguyen Thanh Huong]{Nguyen Thanh Huong}
\address{Thai Nguyen University of Sciences, Thai Nguyen, Viet Nam}
\email{nguyenthanhhuong2806@gmail.com}

\author[Vu Vinh Quang]{Vu Vinh Quang}
\address{Thai Nguyen University of Information and Communication Technology}
\email{vvquang@ictu.edu.vn}


\date{\today}

\begin{document}

\begin{abstract}
This paper presents high-order numerical methods for solving boundary value problems associated with the Lane-Emden equation, which frequently arises in astrophysics and various nonlinear models. A major challenge in studying this equation lies in its singularity at one endpoint. Prior to constructing the numerical methods, we establish the existence and uniqueness of the solution and propose a continuous iterative method. This continuous method is then discretized using the trapezoidal quadrature rule enhanced with correction terms. As a result, we derive three discrete iterative schemes tailored for three specific cases of the Lane-Emden equation.
We rigorously prove that the proposed methods achieve eighth-order accuracy and convergence. A series of numerical experiments is conducted to validate the theoretical findings and demonstrate the accuracy and convergence order of the proposed schemes, which outperform existing methods. These schemes thus provide efficient tools for solving Lane-Emden boundary value problems and can be readily extended to higher-order nonlinear singular models, such as Emden-Fowler equations, which arise in many applications.
\end{abstract}

\maketitle

\section{Introduction}\label{sec1}
In this paper, we consider the following boundary value problem
\begin{equation}\label{eq1}
	\begin{aligned}
		u''(x)+\dfrac{\beta}{x}u'(x)&=f(x,u(x)), \quad 0 < x \leq 1, \\
		u'(0)&=0,\quad u(1)=\alpha,
	\end{aligned}
\end{equation}
where $\alpha$ is a given finite constant, $\beta \geq 1$, and the function $f(x,u)$ is assumed to be sufficiently smooth. 

\par The differential equation in \eqref{eq1} is known as the Lane-Emden equation, which arises in the mathematical modeling of various phenomena in theoretical physics, astrophysics, chemistry, and physiology. Examples include the thermal behavior of spherical gas clouds, isothermal gas spheres, stellar structure theory, thermal explosions, thermionic current theory, steady-state oxygen diffusion in spherical cells, and the distribution of heat sources in the human head. Due to its wide range of applications, the Lane-Emden equation has attracted considerable research attention.

A major difficulty in analyzing this equation lies in its singularity at 
$x=0$. To the best of our knowledge, the existence and uniqueness of solutions to problem \eqref{eq1} were first established in \cite{Chawla87,Pandey97} using the method of lower and upper solutions. However, verifying the required conditions for this method is often challenging, as it involves constructing appropriate lower and upper solutions.

Since those foundational works, most researchers have shifted focus toward numerical and analytical methods for solving nonlinear singular BVPs of the form \eqref{eq1}. Several analytical approaches have been proposed, including the Adomian decomposition method \cite{Hoss,Liao,Singh,Wazwaz}, the variational iteration method \cite{Ghor,Kanth,Wazwaz1}, the differential transformation method \cite{Khan}, and the homotopy analysis method \cite{Singh3}.

In particular, there has been substantial effort devoted to developing high-accuracy numerical methods for solving the Lane-Emden equation. Among these, the finite difference method is one of the most widely used. Early works by Chawla et al. \cite{Chawla,Chawla1} introduced second- and fourth-order difference schemes for singular BVPs of the form \eqref{eq1}. A second-order scheme based on a non-uniform mesh was proposed in \cite{Kumar}, while \cite{Verma} introduced a nonstandard finite difference scheme to address the singularity at $x=0$.

Recently, compact finite difference methods (CFDMs) have gained popularity due to their ability to achieve high-order accuracy with relatively few grid points. Such methods have been applied to nonlinear singular BVPs in, e.g., \cite{Malele,Roul3}, with the scheme in \cite{Malele} attaining an accuracy of $O(h^7)$.

Another widely adopted approach is the use of spline collocation methods. Representative works include \cite{Alam,Caglar,Kada,Kanth1,Roul1,Roul2}, in which various spline and B-spline techniques have been applied. In \cite{Roul1}, the authors developed a standard quintic B-spline collocation method with fourth-order accuracy, and an optimal version achieving sixth-order accuracy. The convergence proofs in these works rely on Green?s function techniques. In the same year, Roul et al. \cite{Roul2} proposed an optimal sixth-order quartic B-spline method for solving \eqref{eq1}. More recently, Roul \cite{Roul4} introduced a fourth-order method based on quintic B-splines for the Lane-Emden equation with Robin boundary conditions:
\begin{equation*}\label{RBC}
	u'(0)=0, \mu u(1) +\sigma u'(1)=B.
\end{equation*}
Additionally, for the same class of boundary conditions, Alam et al. \cite{Alam} developed a fifth-order method in 2021 using a quartic trigonometric B-spline collocation method combined with a quasilinearization approach.

It is worth emphasizing that in most numerical studies of the Lane-Emden equation, the existence and uniqueness of the solution are assumed implicitly, without detailed discussion.

Furthermore, {\it it is notable that no existing numerical methods for the Lane-Emden equation have achieved an order of accuracy higher than seven. This motivates the development of new high-accuracy schemes.}

In this work, we extend the eighth-order iterative framework developed in \cite{Dang2024-1,Dang2024-2,Dang2025-1}, originally applied to third-, fourth-, and second-order ODEs, to construct an eighth-order numerical method for solving problem \eqref{eq1}.

The remainder of this paper is organized as follows. In Section 2, we establish the existence and uniqueness of the solution to problem \eqref{eq1} by reformulating it as an operator equation and proposing a continuous-level iterative scheme. In Section 3, we develop eighth-order discrete iterative methods for the cases 
$\beta =1, \beta =n >1$  and $\beta = r/s,\ r>s$. Section 4 presents several numerical examples to illustrate the eighth-order convergence and to compare our method with existing ones. In Section 5, we briefly consider the extension to Robin boundary conditions, and Section 6 concludes the paper.
\section{Analytical Foundation: Existence and a Continuous Iterative Method} 
Before developing numerical schemes for problem \eqref{eq1}, we first address its qualitative issue and derive a continuous-level iterative method for its solution. Specifically, we reformulate the boundary value problem as a fixed-point problem for an associated nonlinear operator, following the general framework of \cite{A-Long2021}.
\par Let 
$G_0(x,t)$ denote the Green function corresponding to problem \eqref{eq1}. According to \cite{Singh1}, it is defined as follows:
\begin{equation*}
	G_0(x,t)=\begin{cases}
		\ln t, \quad 0 < x \leq t \leq 1,\\
		\ln x, \quad 0 < t \leq x \leq 1
	\end{cases} \quad (\beta = 1),
\end{equation*}
\begin{equation*}
	G_0(x,t)=\dfrac{1}{1-\beta}\begin{cases}
		t^{1-\beta}-1, \quad 0 < x \leq t \leq 1,\\
		x^{1-\beta}-1, \quad 0 < t \leq x \leq 1
	\end{cases} \quad (\beta > 1).
\end{equation*}
Set 
\begin{equation*}
	G_1(x,t)=\dfrac{\partial G_0(x,t)}{\partial x}.
\end{equation*}
Then,
\begin{equation*}
	G_1(x,t)=\begin{cases}
		0, \quad 0 < x \leq t \leq 1,\\
		\dfrac{1}{x}, \quad 0 < t \leq x \leq 1
	\end{cases} \quad (\beta = 1),
\end{equation*}
\begin{equation*}
	G_1(x,t)=\begin{cases}
		0, \quad 0 < x \leq t \leq 1,\\
		x^{-\beta}, \quad 0 < t \leq x \leq 1
	\end{cases} \quad (\beta > 1).
\end{equation*}
From the structure of $G_0$ and $G_1$, one can derive the following bounds:
\begin{equation}\label{eq2}
	\begin{aligned}
		\max_{0\le x\le 1}	\int_0^1 \Big |t^\beta G_0(x,t)\Big|dt &= \dfrac{1}{2(\beta + 1)},\\
		\max_{0\le x\le 1}	\int_0^1 \Big |t^\beta G_1(x,t)\Big|dt &= \dfrac{1}{\beta + 1}.
	\end{aligned}
\end{equation}
Now we define a nonlinear operator $A:C[0,1] \rightarrow C[0,1]$ by
$$(A\varphi)(x)=f(x,u(x)),$$
where $u(x)$ solves the linear problem
\begin{equation}\label{eq3}
	\begin{aligned}
		u''(x)+\dfrac{\beta}{x}u'(x)&=\varphi(x), \quad 0 < x \leq 1, \\
		u'(0)&=0,\quad u(1)=\alpha.
	\end{aligned}
\end{equation}
Analogously as done in \cite{A-Long2021}, it is easy to verify that
the solution of the problem \eqref{eq1} is reduced to a fixed point problem for the operator $A$,  i.e., 
to the solution of the operator equation $A\varphi = \varphi$.
\par Notice that, the solution to the problem \eqref{eq3} can be expressed in integral form as
$$u(x)=u_\varphi (x)=\alpha + \int _0 ^1 t^\beta G_0(x,t)\varphi(t)dt.$$
Therefore,
$$u'(x)=\int _0 ^1 t^\beta G_1(x,t)\varphi(t)dt.$$
From \eqref{eq2} we obtain the estimates
$$\|u\| \leq |\alpha| + \dfrac{1}{2(\beta + 1)}\|\varphi\|,$$
$$\|u'\| \leq \dfrac{1}{\beta + 1}\|\varphi\|,$$
where $\|.\|$ is the uniform norm in the space $C[0, 1]$.
\par Let $M>0$ be arbitrary and define the domain
\begin{equation*} \label{eq4}
	\begin{split}
		\mathcal{D}_M= \Big \{ (x,u) \mid  0\leq x\leq 1, \; |u|\leq |\alpha| + \dfrac{1}{2(\beta + 1)}M \Big \}.
	\end{split}
\end{equation*}
Based on finding the fixed point of the associated operator $A$,  it is easy to prove the following theorem.
\begin{theorem}\label{theo1}
	Suppose that the function $f$ is continuous and there exists constants $M>0$, $L \geq 0$ such that
	\begin{equation*}
		|f(x,u)| \leq M,
	\end{equation*}
	\begin{equation*}
		|f(x,u_2)-f(x,u_1)| \leq L |u_2-u_1|
	\end{equation*}
	\noindent for all $(x,u), (x,u_1), (x,u_2) \in \mathcal{D}_M$ and
	\begin{equation*} 
		q= \dfrac{L}{2(\beta + 1)} < 1.
	\end{equation*}
	Then the problem \eqref{eq1} admits a unique solution $u$, which satisfies the bounds
	\begin{equation*}
		\| u\| \leq |\alpha| + \dfrac{M}{2(\beta + 1)}, \quad \| u'\| \leq \dfrac{M}{\beta + 1}.
	\end{equation*}
\end{theorem}
\par Consider the following continuous-level iterative method for solving the problem \eqref{eq1}, which in essence is the successive iterative method for finding the fixed point of the operator $A$:\\
\noindent 
1. Initialization: Set 
\begin{equation*}
	\varphi _0(x)=f(x,0).
\end{equation*}
2. Iteration: For $k\geq 0,$ compute:
\begin{equation}\label{eq6}
	u_k(x)=\alpha + \int_0^1 t^\beta G_0(x,t)\varphi _k(t)dt.
\end{equation}
3. Update: 
\begin{equation*}
	\begin{aligned}
		\varphi_{k+1}(x)=f(x,u_k(x)).
	\end{aligned}
\end{equation*}
\begin{theorem}\label{theo2}
	Under the assumptions of Theorem \ref{theo1}, the iterative process described above converges uniformly to the solution $u$
	of problem \eqref{eq1}. Moreover, the convergence is geometric with the following estimates:
	\begin{equation*} 
		\begin{split}
			\| u_k-u\| \leq \dfrac{1}{2(\beta + 1)}p_k, \, \| u'_k-u'\| \leq \dfrac{1}{\beta + 1}p_k,\,
		\end{split}
	\end{equation*}
	where
	\begin{equation*}\label{pk}
		p_k=\dfrac{q^k}{1-q}\| \varphi  _1-\varphi  _0\|.
	\end{equation*}
\end{theorem}
\section{Eighth-Order Discretization Schemes via Corrected Quadrature}\label{NumMeth}
To construct eighth-order accurate numerical methods for problem \eqref{eq1}, we discretize the continuous iterative framework established in Section 2. The key idea is to approximate the integral representation of the solution using a corrected trapezoidal rule derived from the Euler-Maclaurin formula, which incorporates endpoint derivatives for enhanced accuracy.

Let $h=(b-a)/N$, where $N$ is a some positive integer, and $t_i=ih, \, i=0,1,...,N.$ For any function $\Phi \in C^{2p}[a,b]$, the Euler-Maclaurin quadrature formula \cite{Hild} has the form
\begin{equation*}
	\int \limits_a^b \Phi (t)dt=T_\Phi(h)-\sum\limits_{l=1}^{p-1}\dfrac{B_{2l}}{(2l)!} \Big(\Phi^{(2l-1)}(b)-\Phi^{(2l-1)}(a) \Big)h^{2l} + O(h^{2p}),
\end{equation*} 
where $B_{2l}$ are the Bernoulli numbers, and
$$T_\Phi(h)=\dfrac{h}{2}(\Phi_0 + \Phi_{N}) + h\sum\limits_{i=1}^{N-1} \Phi_i,$$
with $\Phi_i=\Phi(t_i)$. For the case $p=4$ we obtain a quadrature formula of eighth-order accuracy
\begin{equation}\label{eq5}
	\begin{aligned}
		\int \limits_a^b \Phi (t)dt&=T_\Phi(h)-\dfrac{h^2}{12}  \Big[\Phi'(b)-\Phi'(a) \Big] 
		+ \dfrac{h^4}{720}  \Big[\Phi'''(b)-\Phi'''(a) \Big] 
		\\
		&- \dfrac{h^6}{30240} \Big[\Phi^{(5)}(b)-\Phi^{(5)}(a) \Big] 
		+ O(h^8).
	\end{aligned}
\end{equation}
These correction terms involve the first, third, and fifth derivatives of 
$\Phi$ at the endpoints, which are approximated numerically using finite difference operators.
\par We now present the discretization schemes for three distinct cases of the parameter $\beta$, each requiring tailored treatment of the singularity and kernel structure.\\
1. Integer $\beta =1$,\\
2. Integer $\beta =n >1$, and\\
3. Rational $\beta = r/s,$ where $r,s \in \mathbb N^*$, $r>s$.
\subsection{Case 1: $\beta = 1$}
We apply the formula \eqref{eq5} to compute the integrals in
$$u(x) = \alpha + \int_0^1 tG_0(x,t)\varphi(t)dt , \quad 0 <x \leq 1,$$
and
$$u(0)=\alpha + \int_0^1 t\ln t\varphi(t)dt.$$
Notice that
\begin{equation*}
	\int_0^1 tG_0(x,t)\varphi(t)dt = \ln x \int _0^x t\varphi(t)dt + \int _x^1 t\ln t \varphi(t)dt.
\end{equation*}
Let 
\begin{equation*}
	I_1(x)= \int _0^x t\varphi(t)dt, \quad I_2(x)= \int _x^1 t\ln t \varphi(t)dt, \quad 0 \leq x \leq 1.
\end{equation*}
We have
\begin{equation}\label{eq7}
	\begin{aligned}
		u(x)&=\alpha + \ln x I_1(x) + I_2(x), \quad 0 < x \leq 1,\\
		u(0) &= \alpha + I_2(0).
	\end{aligned}
\end{equation}
Consider $I_2(x)$ with $0 < x \leq 1$. Let $F(x)$ be a continuous and sufficiently smooth function on $[0,1]$ satisfying
\begin{equation*}
	d(xF(x))=x\varphi(x)dx.
\end{equation*}
We have
$$xF(x)=\int _0^x t\varphi(t)dt.$$
Therefore,
$$F(x)=\dfrac{1}{x}\int _0^x t\varphi(t)dt, \quad 0 < x  \leq 1.$$
From the condition that $F$ is right-continuous at $x=0$, we infer that 
\begin{equation*}
	F(0)=\lim_{x\rightarrow 0+} F(x)=\lim_{x\rightarrow 0+}\dfrac{\int _0^x t\varphi(t)dt}{x}=0.
\end{equation*} 
Therefore, 
\begin{equation}\label{eq8}
	F(x)=\begin{cases}\dfrac{1}{x}\int _0^x t\varphi(t)dt, \quad 0 < x  \leq 1,\\
		0, \quad x = 0.
	\end{cases}
\end{equation}
At this point, we have
\begin{equation*}\label{eq9}
	\begin{aligned}
		I_2(x)&=\int _x^1 t \ln t \varphi(t)dt=\int _x^1 \ln t d(tF(t)) \\
		&=t\ln t F(t)\Big|_x^1 - \int_x^1 F(t)dt\\
		&=-x \ln x F(x) - \int_x^1 F(t)dt, \quad 0 < x \leq 1,
	\end{aligned}
\end{equation*}
\begin{equation*}\label{eq10}
	\begin{aligned}
		I_2(0)=\lim _{x \rightarrow 0+}\Big(-x \ln x F(x) - \int_x^1 F(t)dt \Big)
		=-\int _0^1 F(t)dt.
	\end{aligned}
\end{equation*}
Now consider the uniform grid
$$\overline \omega _h = \{x_i=ih, \quad h=1/N, \quad i=0,1,...,N \}.$$
We have
$$I_1(x_0)=I_1(0)=0.$$
At the points $x_i$ $(i=1,2,...,N),$ we compute $I_1(x_i)$ using the formula \eqref{eq5} as follows
\begin{equation*}
	\begin{aligned}
		I_1(x_i)&=\int_0^{x_i} t\varphi(t)dt\\
		&=\dfrac{h}{2}\Big[x_i\varphi(x_i)+2\sum\limits_{j=1}^{i-1}x_j\varphi(x_j) \Big]
		-\dfrac{h^2}{12}\Big[\varphi(x_i)+x_i\varphi'(x_i)-\varphi(x_0) \Big]\\
		&+\dfrac{h^4}{720}\Big[3\varphi''(x_i)+x_i\varphi'''(x_i)-3\varphi''(x_0) \Big] \\
		&-\dfrac{h^6}{30240}\Big[5\varphi^{(4)}(x_i)+x_i\varphi^{(5)}(x_i)-5\varphi^{(4)}(x_0) \Big]+O(h^8).
	\end{aligned}
\end{equation*}
For short we denote $\varphi_i = \varphi(x_i)$. Denote by $D_m^{(n)}\varphi_i$ the difference formula for approximating $n$th derivative $\varphi^{(n)}(x)$ with $m$-order of accuracy at point $x_i$ $(i=0,1,...,N)$. We have 
\begin{equation}\label{eq11}
	I_1(x_i)=L_8(I_1,x_i)\varphi + O(h^8), \quad i=0,1,...,N,
\end{equation} 
where
\begin{equation}\label{eq12}
	L_8(I_1,x_i)\varphi=0, \quad i=0,
\end{equation}
\begin{equation}\label{eq13}
	\begin{aligned}
		&L_8(I_1,x_i)\varphi = \dfrac{h}{2}\Big[x_i\varphi_i+2\sum\limits_{j=1}^{i-1}x_j\varphi_j \Big]
		-\dfrac{h^2}{12}\Big[\varphi_i+x_iD_6^{(1)}\varphi_i-\varphi_0 \Big]\\
		&+\dfrac{h^4}{720}\Big[3D_4^{(2)}\varphi_i+x_iD_4^{(3)}\varphi_i-3D_4^{(2)}\varphi_0 \Big] \\
		&-\dfrac{h^6}{30240}\Big[5D_2^{(4)}\varphi_i+x_iD_2^{(5)}\varphi_i-5D_2^{(4)}\varphi_0 \Big], \quad i=1,2,...,N.
	\end{aligned}
\end{equation} 
Next, we compute $I_2(x_i)$ for $i=0,1,...,N.$ We have
\begin{equation}\label{eq14}
	\begin{cases}
		I_2(x_0)=-\int \limits_0^1 F(t)dt,\\
		I_2(x_i)=-x_i\ln x_iF(x_i)-\int \limits_{x_i}^1 F(t)dt, \quad i=1,2,...,N-1,\\
		I_2(x_N) =0,
	\end{cases}
\end{equation}
where $F$ is defined by \eqref{eq8}, that is,
\begin{equation*}
	F(x_0)=0, \quad F(x_i)=\dfrac{1}{x_i}I_1(x_i), \quad (i=1,2,...,N),
\end{equation*}
$I_1(x_i)$ being defined by \eqref{eq11}-\eqref{eq13}.
\par Using \eqref{eq5}, we obtain
\begin{equation*}
	\begin{aligned}
		\int_{x_i}^1 F(t)dt&=\dfrac{h}{2}\Big[F(x_N)+F(x_i)+2\sum \limits_{j=i+1}^{N-1} F(x_j) \Big]-\dfrac{h^2}{2}\Big[ F'(x_N)-F'(x_i)\Big]\\
		&+\dfrac{h^4}{720}\Big[F'''(x_N)-F'''(x_i)\Big]-\dfrac{h^6}{30240}\Big[F^{(5)}(x_N)-F^{(5)}(x_i)\Big]\\&+O(h^8),\quad i=0,1,...,N-1.
	\end{aligned}
\end{equation*}
Therefore,
\begin{equation}\label{eq15}
	\int_{x_i}^1 F(t)dt=L_8(F,x_i)\varphi + O(h^8),
\end{equation}
where
\begin{equation}\label{eq16}
	\begin{aligned}
		&L_8(F,x_i)\varphi=\dfrac{h}{2}\Big[F_N+F_i+2\sum \limits_{j=i+1}^{N-1} F_j \Big]-\dfrac{h^2}{2}\Big[ D_6^{(1)}F_N-D_6^{(1)}F_i\Big]\\
		&+\dfrac{h^4}{720}\Big[D_4^{(3)}F_N-D_4^{(3)}F_i\Big]-\dfrac{h^6}{30240}\Big[D_2^{(5)}F_N-D_2^{(5)}F_i\Big],\quad i=0,1,...,N-1.
	\end{aligned}
\end{equation}
In the above formula $F_i = F(x_i)$ and $D_m^{(n)}F_i$ is the difference formula for approximating $n$th derivative $F^{(n)}(x)$ with $m$-order accuracy at the point $x_i$.
From \eqref{eq14}-\eqref{eq16} we obtain 
\begin{equation}\label{eq17}
	I_2(x_i)=L_8(I_2,x_i)\varphi + O(h^8), \quad i = 0,1,...,N,
\end{equation}
where
\begin{equation}\label{eq18}
	\begin{cases}
		L_8(I_2,x_i)\varphi=-L_8(F,x_0)\varphi, \quad i = 0,\\
		L_8(I_2,x_i)\varphi=-x_i\ln x_i F_i-L_8(F,x_i)\varphi, \quad i=1,2,...,N-1,\\
		L_8(I_2,x_i)\varphi=0, \quad i=N.
	\end{cases}
\end{equation}
Denote by $\Phi_k (x), U_k (x)$ the grid functions, which are defined
on the grid $\overline \omega_h$ and approximate the functions $\varphi_k (x), u_k (x)$ on this grid, respectively.
Combining \eqref{eq7}, \eqref{eq11}-\eqref{eq13}, \eqref{eq17} and \eqref{eq18}, we obtain the discrete iterative method in the case $\beta=1$ as follows:\\

\noindent {\bf Method 1:}\\
1. Initialization: Given 
\begin{equation*}\label{eq19}
	\Phi_0(x_i)=f(x_i,0), \quad i=0,1,...,N.
\end{equation*}
2. Iteration: Knowing $\Phi_k(x_i), \, k=0,1,..., \, i =0,1,...,N$, compute approximately the definite integral \eqref{eq6} by the Euler-Maclaurin formula
\begin{equation*}\label{eq20}
	\begin{aligned}
		&U_k(x_i)= \alpha + L_8(I_2,x_i)\Phi_k + \ln x_i L_8(I_1,x_i)\Phi_k, \quad i =1,2,...,N\\
		&U_k(x_0)=\alpha + L_8(I_2,x_0)\Phi_k.
	\end{aligned}
\end{equation*}
3. Update:
\begin{equation*}\label{eq21}
	\Phi_{k+1}(x_i)=f(x_i,U_k(x_i)), \quad i=0,1,...,N.
\end{equation*}
\subsection{Case 2: $\beta = n >1$, integer} \label{sec 2}
We have
\begin{equation*}\label{eq22}
	\begin{aligned}
		&u(x) = \alpha + \int_0^1 t^\beta G_0(x,t)\varphi(t)dt \\
		&=\alpha + \dfrac{1}{1-n}\Big[(x^{1-n}-1)\int_0^x t^n \varphi(t)dt - \int_x^1 t^n\varphi(t)dt + \int_x^1 t\varphi(t)dt \Big], \quad 0 <x \leq 1,
	\end{aligned}
\end{equation*}
and
\begin{equation*}\label{eq23}
	u(0)=\alpha + \dfrac{1}{1-n}\Big[\int_0^1 t\varphi(t)dt-\int_0^1 t^n\varphi(t)dt\Big].
\end{equation*}
Denote 
\begin{equation*}\label{eq24}
	A_1(x)= \int _0^x t^n\varphi(t)dt, \quad A_2(x)= \int _x^1 t^n \varphi(t)dt, \quad A_3(x)= \int _x^1 t \varphi(t)dt, \quad 0 \leq x \leq 1.
\end{equation*}
We obtain
\begin{equation}\label{eq28}
	\begin{aligned}
		u(x)&=\alpha + \dfrac{1}{1-n}\Big[(x^{1-n}-1) A_1(x) - A_2(x)+A_3(x)\Big], \quad 0 < x \leq 1,\\
		u(0) &= \alpha + \dfrac{1}{1-n}[A_3(0)-A_2(0)].
	\end{aligned}
\end{equation}
Similar to the previous subsection, using the Euler-Maclaurin formula to compute the integrals $A_1(x), A_2(x), A_3(x)$ with accuracy of order 8, we obtain 
\begin{equation}\label{eq31}
	A_3(x_i)=\int_{x_i}^1 t\varphi(t)dt = L_8(A_3,x_i)\varphi + O(h^8),
\end{equation}
where
\begin{equation}\label{eq32}
	L_8(A_3,x_i)\varphi = 0, \quad i = N,
\end{equation}
\begin{equation}\label{eq33}
	\begin{aligned}
		&L_8(A_3,x_i)\varphi = \dfrac{h}{2}\Big[x_i\varphi_i + x_N\varphi_N +2\sum \limits_{j=i+1}^{N-1}x_j\varphi_j \Big]\\
		&-\dfrac{h^2}{12}\Big[x_N D_6^{(1)}\varphi_N +\varphi_N- x_iD_6^{(1)}\varphi_i -\varphi_i \Big]\\
		&+\dfrac{h^4}{720}\Big[3 D_4^{(2)}\varphi_N + x_N D_4^{(3)}\varphi_N -3 D_4^{(2)}\varphi_i - x_i D_4^{(3)}\varphi_i \Big]\\
		&-\dfrac{h^6}{30240}\Big[x_N D_2^{(5)}\varphi_N + 5 D_2^{(4)}\varphi_N -x_i D_2^{(5)}\varphi_i -5D_2^{(4)}\varphi_i \Big], \quad i =0,1,..., N-1.
	\end{aligned}
\end{equation}
Consider
$$A_1(x_i)=\int_0^{x_i} t^n\varphi(t)dt.$$
Denote
$$d_i^p=(t^n)^{(p)}\Big|_{t=x_i}=\begin{cases}
	x_i^n, \quad p = 0,\\
	0, \quad p > n,\\
	n(n-1)...(n-p+1)x_i^{n-p}, \quad 0 <p \leq n.
\end{cases}$$
We have
\begin{equation}\label{eq34}
	A_1(x_i)= L_8(A_1,x_i)\varphi + O(h^8),
\end{equation}
where
\begin{equation}\label{eq35}
	L_8(A_1,x_i)\varphi = 0, \quad i = 0,
\end{equation}
\begin{equation}\label{eq36}
	\begin{aligned}
		&L_8(A_1,x_i)\varphi = \dfrac{h}{2}\Big[x_i^n\varphi_i + x_0^n\varphi_0 +2\sum \limits_{j=1}^{i-1}x_j^n\varphi_j \Big]\\
		&-\dfrac{h^2}{12}\Big[d_i^0 D_6^{(1)}\varphi_i +d_i^1\varphi_i- d_0^0D_6^{(1)}\varphi_0 -d_0^1\varphi_0 \Big]\\
		&+\dfrac{h^4}{720}\Big[\sum_{m=0}^3 C^3_m\Big(d_i^m D_4^{(3-m)}\varphi_i  -d_0^m D_4^{(3-m)}\varphi_0 \Big) \Big]\\
		&-\dfrac{h^6}{30240}\Big[\sum_{m=0}^5 C^5_m\Big(d_i^m D_2^{(5-m)}\varphi_i  -d_0^m D_2^{(5-m)}\varphi_0 \Big) \Big], \quad i =1,2,..., N.
	\end{aligned}
\end{equation}
Similarly, we obtain
\begin{equation}\label{eq37}
	A_2(x_i)= \int_{x_i}^1 t^n\varphi(t)dt=L_8(A_2,x_i)\varphi + O(h^8),
\end{equation}
where
\begin{equation}\label{eq38}
	L_8(A_2,x_i)\varphi = 0, \quad i = N,
\end{equation}
\begin{equation}\label{eq39}
	\begin{aligned}
		&L_8(A_2,x_i)\varphi = \dfrac{h}{2}\Big[x_N^n\varphi_N + x_i^n\varphi_i +2\sum \limits_{j=i+1}^{N-1}x_j^n\varphi_j \Big]\\
		&-\dfrac{h^2}{12}\Big[d_N^0 D_6^{(1)}\varphi_N +d_N^1\varphi_N- d_i^0D_6^{(1)}\varphi_i -d_i^1\varphi_i \Big]\\
		&+\dfrac{h^4}{720}\Big[\sum_{m=0}^3 C^3_m\Big(d_N^m D_4^{(3-m)}\varphi_N  -d_i^m D_4^{(3-m)}\varphi_i \Big) \Big]\\
		&-\dfrac{h^6}{30240}\Big[\sum_{m=0}^5 C^5_m\Big(d_N^m D_2^{(5-m)}\varphi_N  -d_i^m D_2^{(5-m)}\varphi_i \Big) \Big], \quad i =0,1,..., N-1.
	\end{aligned}
\end{equation}
Combining \eqref{eq28}-\eqref{eq39}, we obtain the discrete iterative method in the case $\beta=n \in \mathbb N^*, \, n>1$ as follows:\\

\noindent {\bf Method 2:}\\
1. Initialization: Given 
\begin{equation*}\label{eq40}
	\Phi_0(x_i)=f(x_i,0), \quad i=0,1,...,N.
\end{equation*}
2. Iteration: Knowing $\Phi_k(x_i), \, k=0,1,..., \, i =0,1,...,N$, compute approximately the definite integral \eqref{eq6} by the Euler-Maclaurin formula
\begin{equation*}\label{eq41}
	\begin{aligned}
		&U_k(x_i)= \alpha + \dfrac{1}{1-n}\Big[(x_i^{1-n}-1)L_8(A_1,x_i)\Phi_k - L_8(A_2,x_i)\Phi_k + L_8(A_3,x_i)\varphi_k \Big],\\
		&\text {for} \quad i =1,2,...,N,\\
		&U_k(x_0)=\alpha + \dfrac{1}{1-n}\Big[L_8(A_3,x_0)\Phi_k - L_8(A_2,x_0)\varphi_k \Big].
	\end{aligned}
\end{equation*}
3. Update:
\begin{equation*}\label{eq42}
	\Phi_{k+1}(x_i)=f(x_i,U_k(x_i)), \quad i=0,1,...,N.
\end{equation*}
\subsection{Case 3: $\beta = \frac{r}{s},\, (r,s \in \mathbb N^*, \, r>s$), rational}
In this case, the problem \eqref{eq1} becomes
\begin{equation}\label{eq1A}
	\begin{aligned}
		u''(x)+\dfrac{r/s}{x}u'(x)&=f(x,u(x)), \quad 0 < x \leq 1, \\
		u'(0)&=0,\quad u(1)=\alpha,
	\end{aligned}
\end{equation}
Let $y=x^{1/s}$ and $v(y)=u(y^s)=u(x)$. We deduce that $0<y\leq 1$ because $0<x\leq 1.$\\
We have 
$$y'_x=\dfrac{1}{s}y^{1-s},$$
so
\begin{align}
	u'_x&=v'_y.y'_x=\dfrac{1}{s}y^{1-s}v'_y,\nonumber\\
	u''_{xx}&=\Big(\dfrac{1}{s}y^{1-s}v'_y\Big)'_yy'_x 
	=\dfrac{1}{s^2}y^{2-2s}v''_{yy}+\dfrac{1-s}{s^2}y^{1-2s}v'_y. \nonumber
\end{align}
Then, the problem \eqref{eq1A} (which is also problem \eqref{eq1}) with the variable $x$ is transformed into a new problem with the variable $y$ as follows:
\begin{equation}\label{eq1B}
	\begin{aligned}
		v''(y)+\dfrac{r-s+1}{y}v'(y)&=s^2y^{2s-2}f(y^s,v(y)), \quad 0 < y \leq 1, \\
		v'(0)&=0,\quad v(1)=\alpha.
	\end{aligned}
\end{equation}
\begin{remark}
	Notice that $r-s+1=n \in \mathbb N^*, \, n>1.$
	Therefore, with the variable substitution $y=x^{1/s},$ the problem \eqref{eq1A} is reduced to the problem \eqref{eq1B}, which is of type \eqref{eq1} with the parameter $\beta=n=r-s+1 \in \mathbb N^*, \, n>1$. So, we can apply the results obtained in subsection \ref{sec 2} to this problem.
\end{remark}
For the problem \eqref{eq1A} consider the grid
$$\overline \omega_{h}=\{0=x_0<x_1<...<x_n=1\},$$
and  for the problem \eqref{eq1B} the grid
$$\omega_{h_v}=\{y_i=ih_v, \, h_v=1/N, \, (i=0,1,...,N)\},$$
which satisfies $$x_i=(y_i)^{s}, \quad i=0,1,...,N.$$
We obtain the numerical algorithm for solving the problem \eqref{eq1A} as follows:\\

\noindent {\bf Method 3:}\\ 
Step 1: Using Method 2 to find the solution of problem \eqref{eq1B} on the grid $\omega_{h_v}$. Let $V_k$ denote the numerical solution vector of the problem \eqref{eq1B} at the $k$-th iteration.\\
Step 2: Determine the numerical solution of the problem \eqref{eq1A} on the grid 
$\overline \omega_{h}$ from the numerical solution of the problem \eqref{eq1B}  obtained in Step 1
$$U_k(x_i) = V_k(y_i), \quad i = 0,1,...,N.$$
Now we study the convergence of the above iterative methods.
\begin{proposition}\label{pro1}
	Under the assumption that the function $f(t,u)$ is sufficiently smooth, for Method 1, 2, 3 we have for $k=0,1,...$
	$$\|U_k-u_k\|=O(h^8),$$
	where $\|.\|=\|.\|_{\overline \omega_h}$ is the max-norm of function on the grid $\overline \omega_h$.
\end{proposition}
\noindent {\bf Proof.}
The proposition can be proved by induction in a similar way as Proposition 4 in \cite{A-Long2021}.\\
Now combining the above proposition and Theorem \ref{theo2} we obtain the following theorem.
\begin{theorem}\label{theo3}
	Let the assumptions of Proposition \ref{pro1} and Theorem \ref{theo2} hold. Then, for the numerical solution $U_k$ of problem \eqref{eq1} computed by the discrete method Method $r$, where $r \in \{1,2,3\}$, the following error estimates hold:
	$$\|U_k-u\|\leq \dfrac{1}{2(\beta+1)} p_k +O(h^8).$$
\end{theorem}
\section{Numerical Experiments and Performance Evaluation}
This section focuses on evaluating the computational performance of the proposed eighth-order schemes for solving the Lane-Emden equation across a variety of parameter regimes. To this end, we present a series of benchmark examples, including both problems with known exact solutions and those for which the exact solution is not available. This dual approach enables a comprehensive assessment of the accuracy and convergence properties of the methods.

In the numerical tables that follow, the reported error $E$ is defined as $E=\| U_k-u\|$, when the exact solution $u$ is known. For problems without an exact solution, we employ the double-mesh technique and compute
$$E=        \max _{0\le i\le N}|U_N(x_i)-U_{2N(x_{2i})}|.$$
The empirical convergence order is calculated by
$$Order=\log_2 \dfrac{\| U_k^{N/2}-u\|}{\| U_k^N-u\|},$$ 
where the superscripts $N/2$ and $N$ indicate the number of grid points used in the respective approximations.\\
The iteration process is continued until the stopping criterion
$$\| \Phi_{k+1}-\Phi_k\| \leq TOL =10^{-22}$$
is satisfied.\par
All computations are carried out in MATLAB (version 7.5.0) on a Sony Vaio laptop equipped with an Intel(R) Core(TM) i3-3120M CPU @ 2.50GHz and 4.00 GB of RAM.\\
\noindent {\bf Example 1.} (Example 1 in  \cite{Roul1}, Example 3 in \cite{Roul2})\\
Consider the problem
\begin{equation*}\label{exam2}
	\left\{
	\begin{array}{lll}
		u''(x)+\dfrac{1}{x}u'(x)=e^{u(x)}, \quad 0 < x \leq 1,\\
		u'(0)=0, \quad u(1)= 0.
	\end{array}\right.
\end{equation*}
The exact solution of the problem is
$u(x)=2\ln \Big(\dfrac{d+1}{dx^2+1} \Big),$ where $d=-5+2\sqrt 6.$
For this example, we have $\beta = 1, \, \alpha = 0$,
the domain 
\begin{equation*}
	\begin{aligned}
		\mathcal{D}_M&=\Big\{(x,u) \; | \; 0 \leq x \leq 1, \; |u|\leq |\alpha|+\frac{M}{2(\beta + 1)} \Big\}\\
		&=\Big\{(x,u) \; | \; 0 \leq x \leq 1, \; |u|\leq \frac{M}{4} \Big\}. 
	\end{aligned}
\end{equation*}
In the domain $\mathcal{D}_M$, the right-hand side function $|f(x,u)|=|e^u| \leq e^{M/4}$. Therefore, we can take $M=4$ so that $|f(x,u)| \le M$ in $\mathcal{ D}_M$. In the domain $\mathcal{D}_4$, we have $|f'_u|=|e^u|\leq e.$ Hence, we can choose the Lipschitz coefficient $L=e$. In this case $q=\frac{L}{2(\beta + 1)}=\frac{e}{4} <1$. Thus, all the conditions of Theorem \ref{theo1} are satisfied. Hence, the problem under consideration has a unique solution.

To compute the approximate solution of the problem, we use Method 1.
The number of iterations performed for achieving the given $TOL$ is 27.
The results of the computation, in comparison with the optimal quintic B-spline collocation method (OQM) \cite{Roul1} and the optimal quartic B-spline collocation method (OQBCM) \cite{Roul2}, are presented in Table \ref{Tab1}.\par
\begin{table}[!ht]
	\centering
	\caption[smallcaption]{The errors  and order of convergence of Method 1 for Example 1 in comparison with OQM and OQBCM }
	\label{Tab1}
	\begin{tabular}{ ccccccc} 
		\hline
		$N$& $E$    &$ Order$& OQM     &$ Order$ & OQBCM&$ Order$  \\
		\hline
		8&3.8713e-10 &      &8.5381e-10&      &3.5926e-09&\\
		16&9.1226e-13&8.7292&2.1910e-11&5.2842&5.4654e-11&6.038\\
		32&1.4451e-15&9.3021&3.9240e-13&5.8031&7.1101e-13&6.264\\
		64&8.8928e-18&7.3443&6.3838e-15&5.9417&1.1102e-14&6.000\\
		\hline
	\end{tabular}
\end{table}
\noindent {\bf Example 2.} ( Example 4.3 in \cite{Singh1}, Example 2 in \cite{Alam} ) \\
Consider the problem
\begin{equation}\label{eqexam2}
	\left\{
	\begin{array}{lll}
		u''(x)+\dfrac{\beta}{x}u'(x)=-\mu e^{u(x)}, \quad 0 < x \leq 1,\\
		u'(0)=0, \quad u(1)= 0.
	\end{array}\right.
\end{equation}
It is possible to verify that for $\mu \le 1.5$, all the conditions of Theorem \ref{theo1} are satisfied, so the problem has a unique solution.
In the case $\beta =1$,  $\mu =0.5$, the exact solution is not known. The results of computation  by Method 1 are reported in Table \ref{Tab2a}. 
\begin{table}[!ht]
	\centering
	\caption[smallcaption]{The errors and order of convergence of Method 1 for Example 2 in the case $\beta = 1, \mu = 0.5$ }
	\label{Tab2a}
	\begin{tabular}{ cccc} 
		\hline
		$N$ & $k$& $E$ & $ Order$ \\ 
		\hline
		8&22&6.6589e-12&          \\ 
		16&22&7.2623e-15&9.8406\\
		32&22&1.4189e-17&8.9995\\
		64&22&3.8579e-20&8.5227\\
		\hline
	\end{tabular}
\end{table}
In the case $\beta =1, \mu =1$, the exact solution of the problem is $u(x)=2\log\frac{S+1}{Sx^2+1}$, where $S=3-2\sqrt{2}$ (see \cite{Alam}). The errors and order of convergence of Method 1 for this case, in comparison with those of the method proposed by Alam et al. \cite{Alam} are presented in Table \ref{Tab2b}.

\begin{table}[!ht]
	\centering
	\caption[smallcaption]{The errors and order of convergence of Method 1 for Example 2 in the case $\beta = 1, \mu = 1$ in comparison with \cite{Alam} } 
	\label{Tab2b}
	\begin{tabular}{ccccc}
		\hline
		$N$   & $E$ & $ Order$ & $E$ in \cite{Alam} & $ Order$ \\ \hline
		8 & 6.0703e-10 & ~ & ~ & ~ \\ 
		16 & 9.9952e-13 & 9.2463 & 1.34662e-08 & ~ \\ 
		32 & 1.0878e-15 & 9.8437 & 7.86881e-10 & 4.0971 \\ 
		64 & 1.8019e-18 & 9.2377 & 4.76979e-11 & 4.0442  \\ \hline
	\end{tabular}
\end{table}

We also solve the problem \eqref{eqexam2} for the case $\beta =2$ for $\mu =0.5, 1$ by Method 2. The errors and convergence rate are given in Table \ref{Tab2c}.

\begin{table}[!ht]
	\centering
	\caption[smallcaption]{The errors and order of convergence of Method 2 for Example 2 in the case $\beta = 2$ }
	\label{Tab2c}
	\begin{tabular}{ ccccccc} 
		\hline
		$N$ & $k$ & $E\, (\mu=0.5)$ & $Order$ & $k$ & $E\, (\mu=1)$& $Order$ \\
		\hline
		8&18&2.6804e-13&       &24&1.7535e-11& \\
		16&18&9.5929e-16&8.1263&24&3.6259e-14&8.9177\\
		32&18&3.6359e-18&8.0435&24&1.4861e-16&7.9307\\
		64&18&1.7545e-20&7.6951&24&5.6253e-19&8.0454\\
		\hline
	\end{tabular}
\end{table}
\noindent {\bf Example 3.} (Example 2 in \cite{Roul1}, Example 4 in \cite{Roul2})\\
Consider the problem describing the equilibrium of the isothermal gas sphere
\begin{equation*}
	\begin{cases}
		u''(x)+\dfrac{2}{x}u'(x)=-u^{(5)}(x), \quad 0 < x \leq 1,\\
		u'(0)=0, \quad u(1)= \sqrt {3/4}.
	\end{cases}
\end{equation*}
The exact solution of the problem is $u(x)=\sqrt{\dfrac{3}{3+x^2}}$. For this example, $\beta = 2$.

We apply Method 2 described in Section \ref{NumMeth} to solve the problem.
The errors and order of convergence of the method for the example compared with those of OQM in \cite{Roul1} and  OQBCM \cite{Roul2} are given in Table \ref{Tab3}.

\begin{table}[!ht]
	\centering
	\caption[smallcaption]{The errors and order of convergence of Method 2 for Example 3 in comparison with OQM \cite{Roul1} and OQBCM \cite{Roul2}}
	\label{Tab3}
	\begin{tabular}{ccccccc}
		\hline
		$N$ &$E$ & $ Order$  & OQM & $ Order$ & OQBCM & $ Order$ \\ \hline
		8 & 2.5566e-09 & ~ & 4.1099e-08 & ~ & 5.3361e-08 & ~ \\ 
		16 & 5.2076e-12 & 8.9394 & 6.9891e-10 & 5.8778 & 1.5948e-09 & 5.064 \\ 
		32 & 1.0483e-14 & 8.9564 & 1.0676e-11 & 6.0237 & 2.6749e-11 & 5.897 \\ 
		64 & 3.7897e-17 & 8.1118 & 1.6513e-13 & 6.0146 & 4.0890e-13 & 6.031 \\ \hline
	\end{tabular}
\end{table}

\begin{figure}[!ht]
	\begin{minipage}[b] {1.0\textwidth}
		\centering
		\includegraphics[height=5.55cm,width=9cm]{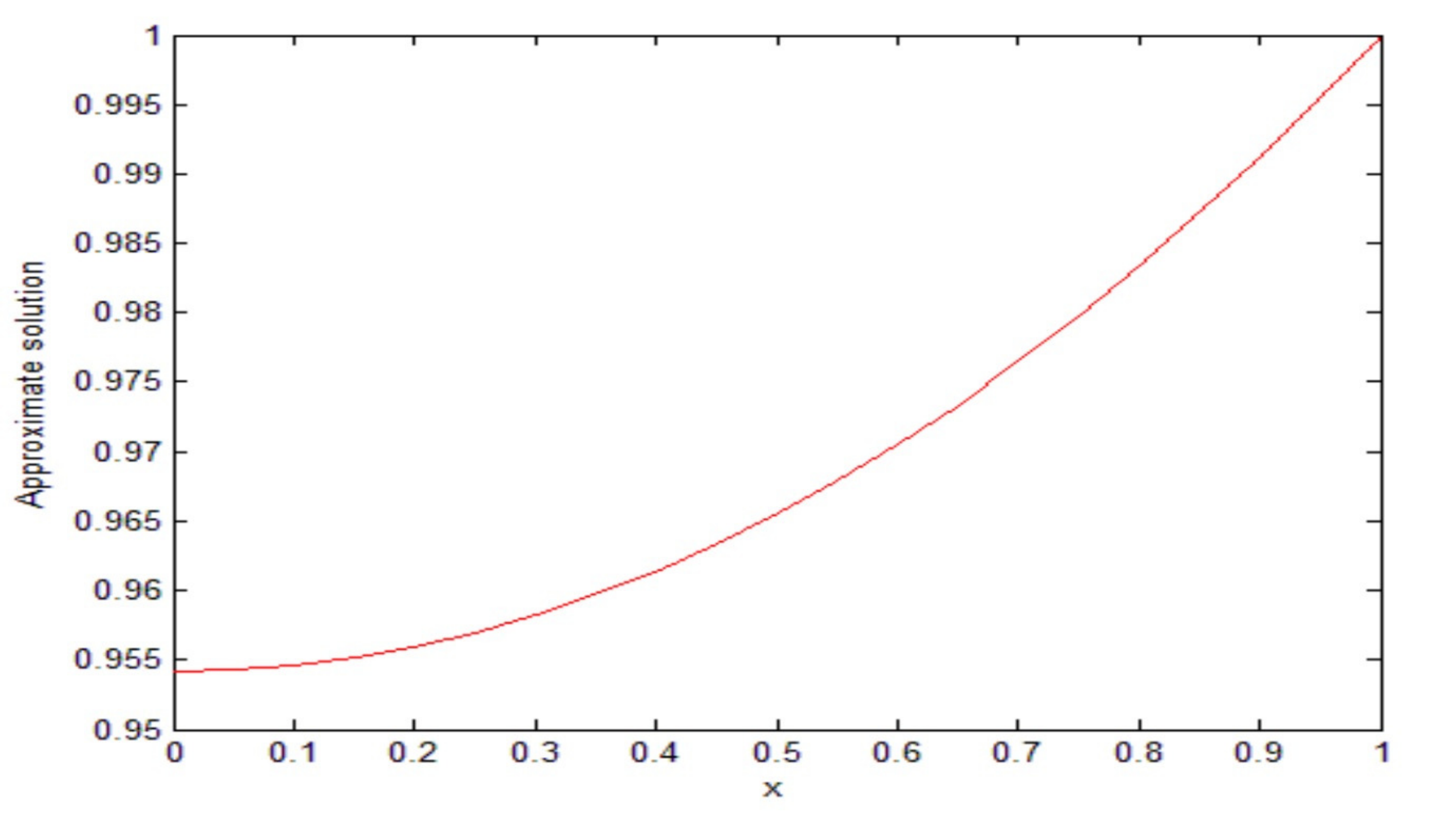}
		\caption{The graph of the approximate solution in Example 4 obtained by Method 2.}\label{figA5}
	\end{minipage}
	\vfill
	\begin{minipage}[b]{1.0\textwidth}
		\centering
		\includegraphics[height=6cm,width=9cm]{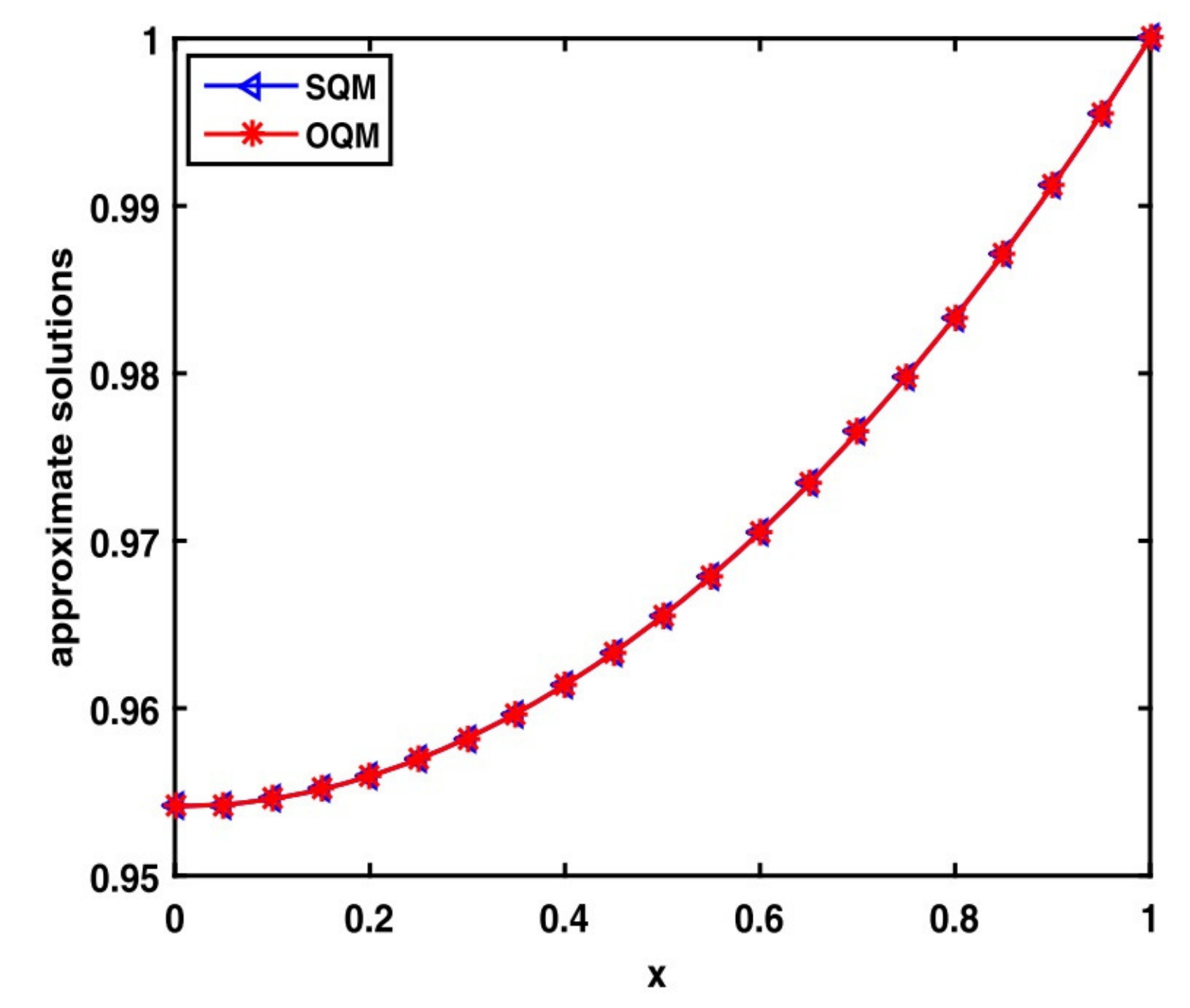}
		\caption{The graph of the approximate solution in Example 4 obtained in \cite{Roul1}}\label{figR5}
	\end{minipage}
\end{figure}

\noindent {\bf Example 4.} (Example 4 in \cite{Roul1})
Consider the nonlinear problem describing the distribution of radial stress on a rotationally symmetric shallow membrane cap \cite{Dickey}
\begin{equation*}\label{exam5}
	\left\{
	\begin{array}{lll}
		(x^3u'(x))'=x^3\Big(\dfrac{1}{2}-\dfrac{1}{8u^2(x)} \Big), \, 0 < x \leq 1,\\
		u'(0)=0, \quad u(1)= 1.
	\end{array}\right.
\end{equation*}
This problem corresponds to \eqref{eq1} with $\beta =3, \alpha =1$ and $f(x,u)=\dfrac{1}{2}-\dfrac{1}{8u^2(x)}$.
The exact solution of the problem is not known. The errors and order of convergence computed by Method 2 are given in Table \ref{Tab5}. The graph of the approximate solution by Method 2 for $N=20$ and the graphic taken from \cite{Roul1} are depicted in Figures \ref{figA5} and \ref{figR5}, respectively. 
\begin{table}[!ht]
	\centering
	\caption[smallcaption]{The errors and order of convergence of Method 2 for Example 4}
	\label{Tab5}
	\begin{tabular}{ cccc} 
		\hline
		$N$ & $k$& $E$ & $ Order$ \\ 
		\hline
		8&12&1.1595e-13&          \\ 
		16&12&4.5242e-16&8.0016\\
		32&12&1.7560e-18&8.0092\\
		64&12&6.8507e-21&8.0018\\
		\hline
	\end{tabular}
\end{table}

\begin{figure}[ht]
	\begin{center}
		\includegraphics[height=6cm,width=9cm]{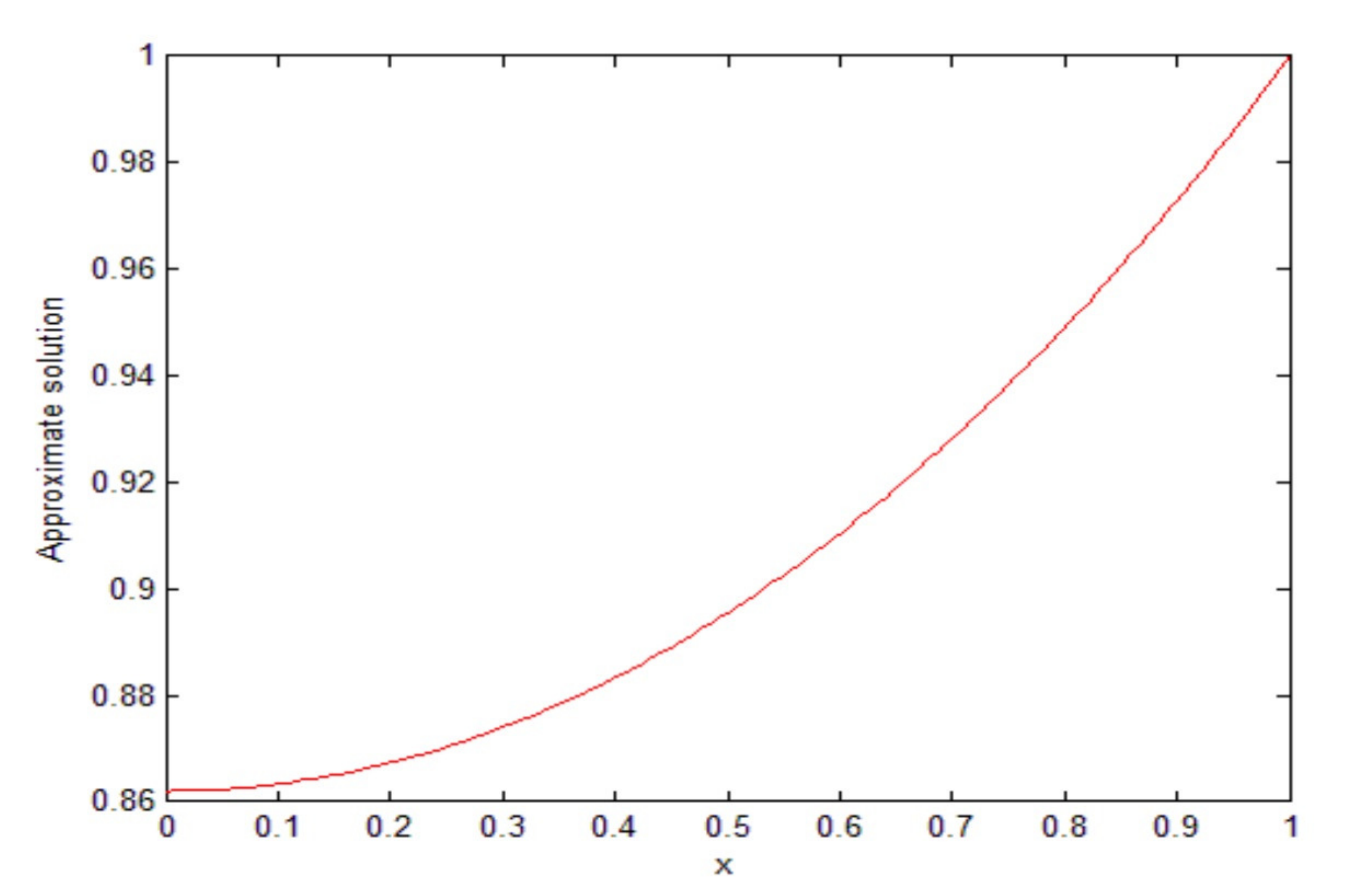}
		\caption{The graph of the approximate solution in Example 5}
		\label{figA6}
	\end{center}
\end{figure}

\noindent {\bf Example 5.} (Example 4.3 in \cite{Singh2})
Consider the problem which arises in the study of electroactive polymer film \cite{Lions}
\begin{equation*}
	\left\{
	\begin{array}{lll}
		(x^2u'(x))'=\dfrac{\mu u(x)}{1+\lambda u(x)}, \quad 0 < x \leq 1,\\
		u'(0)=0, \quad u(1)= 1.
	\end{array}\right.
\end{equation*}
For this problem, $\beta=2$ and the exact solution is unknown in literature. We apply the proposed Method 2 to approximate solution of the problem for the case $\mu = 1, \lambda =0.1$.
The results of errors and convergence rate  in comparison with the CFDM in \cite{Singh2} are given in Table \ref{Tab6a}. The graph of the approximate solution is depicted in Figure \ref{figA6}.

\begin{table}[!ht]
	\centering
	\caption[smallcaption]{The errors and order of convergence for Example 5}
	\label{Tab6a}
	\begin{tabular}{ccccc}
		\hline
		$N$  & $E$ & $ Order$ & $E$ in \cite{Singh2} & $ Order$ \\ \hline
		8 & 1.0617e-13 & ~ & ~ & ~ \\ 
		16 & 4.1470e-16 & 8.0001 & 1.74e-08 & ~ \\ 
		32 & 1.6457e-18 & 7.9772 & 1.20e-09 & 3.85678 \\ 
		64 & 8.6973e-21 & 7.5639 & 4.52e-12 & 4.00919 \\ \hline
	\end{tabular}
\end{table}

\medbreak
\noindent {\bf Remark:} Above, we consider some examples of boundary value problems for Lane-Emden equation \eqref{eq1} with integer $\beta \ge 1$. The numerical results show the eighth order convergence of Method 1 for $\beta =1$ and Method 2 for $\beta >1$. Below, we consider an example of Lane-Emden equation with fractional $\beta >1$.\\

\noindent {\bf Example 6.} 
Consider the problem
\begin{equation}\label{exam8}
	\begin{cases}
		(x^{3/2}u'(x))'=x^{3/2}(f(x)+g(x)e^{u(x)}), \quad 0 < x \leq 1,\\
		u'(0)=0, \quad u(1)= 0,
	\end{cases}
\end{equation}
where
$$f(x)=\dfrac{1}{(x^2+1)^{3/2}}+\dfrac{3/2}{\sqrt{x^2+1}}-1, \quad g(x)=e^{-\sqrt{x^2+1}}.$$
The exact solution of the problem is $u(x)=\sqrt {x^2 + 1}$. \\
Let $t=x^{1/2}$, then the problem \eqref{exam8} becomes
\begin{equation}\label{exam8'}
	\begin{cases}
		v''(t)+\dfrac{2}{t}v'(t)=4t^2(F(t)+G(t)e^{v(t)}), \quad 0 < t \leq 1,\\
		v'(0)=0, \quad v(1)=0,
	\end{cases}
\end{equation}
where $$F(t)=\dfrac{1}{(t^4+1)^{3/2}}+\dfrac{3/2}{\sqrt{t^4+1}}-1, \quad G(t)=e^{-\sqrt{t^4+1}}.$$
Clearly, the problem \eqref{exam8} is reduced to the problem \eqref{exam8'}, which is of the type \eqref{eq1} for $\beta = 2$. Therefore, to solve problem \eqref{exam8}, we apply Method 3 described in Section~\ref{NumMeth}, where Method 2 is employed in the first step. The errors and order of convergence of the method are presented in Table \ref{Tab14}.\\

\begin{table}[!ht]
	\centering
	\caption[smallcaption]{The errors and order of convergence of Method 3 for Example 6  }
	\label{Tab14}
	\begin{tabular}{ cccc} 
		\hline
		$N$ & $k$& $E$ & $ Order$ \\ 
		\hline
		8&21&3.3466e-06&          \\ 
		16&21&1.6565e-09&10.9803\\
		32&21&1.8807e-12&9.7827\\
		64&21&1.7764e-15&10.0481\\
		\hline
	\end{tabular}
\end{table}
\section{Extension to Robin-Type Lane-Emden Problems}
Now consider the problem
\begin{equation}\label{eq1R}
	\begin{aligned}
		u''(x)+\dfrac{\beta}{x}u'(x)&=f(x,u(x)), \quad 0 < x \leq 1, \\
		u'(0)&=0,\quad \mu u(1) +\sigma u'(1)=\alpha ,
	\end{aligned}
\end{equation}
where $\mu >0, \sigma \geq 0$.
For this problem, the Green function has the form (see \cite{Singh})
\begin{equation*}
	G_0(x,t)=\begin{cases}
		\ln t -\frac{\sigma}{\mu}, \quad 0 < x \leq t \leq 1,\\
		\ln x -\frac{\sigma}{\mu}, \quad 0 < t \leq x \leq 1
	\end{cases} \quad (\beta = 1),
\end{equation*}
and
\begin{equation*}
	G_0(x,t)=\begin{cases}
		\frac{t^{1-\beta}-1}{1-\beta} -\frac{\sigma}{\mu}, \quad 0 < x \leq t \leq 1,\\
		\frac{x^{1-\beta}-1}{1-\beta} -\frac{\sigma}{\mu}, \quad 0 < t \leq x \leq 1
	\end{cases} \quad (\beta > 1).
\end{equation*}
By applying the same approach to the problem \eqref{eq1}, we can construct eighth-order numerical method for the problem \eqref{eq1R}. Below is an example for illustrating the eighth-order accuracy of the method for solving Lane-Emden equation with Robin boundary conditions.\\

\noindent {\bf Example 7.} (Example 7 in \cite{Alam}, Example 6.6 in \cite{Malele})
Consider the problem with Robin conditions
\begin{equation*}\label{exam10}
	\left\{
	\begin{array}{lll}
		u''(x)+\dfrac{2}{x}u'(x)=-e^{-u(x)}, \quad 0 < x \leq 1,\\
		u'(0)=0, \quad 2u(1)+u'(1)= 0
	\end{array}\right.
\end{equation*}
which arises in the study of thermal distribution in the human head \cite{Duggan}.
The exact solution of the problem is not known. We solved numerically the problem by Method 2 with some adaptation.
Table \ref{Tab10} presents the errors (computed by double mesh principle) and the order of convergence of the method. 

\begin{table}[!ht]
	\centering
	\caption[smallcaption]{The errors and order of convergence in  Example 7  }
	\label{Tab10}
	\begin{tabular}{ccc}
		\hline
		$N$  & $E$ & $ Order$   \\ \hline
		8 & 6.3628e-11 & ~ \\ 
		16 & 1.3235e-13 & 8.9092 \\ 
		32 & 2.1263e-16 & 9.2818 \\ 
		64 & 5.7366e-19 & 8.5339 \\ \hline
	\end{tabular}
\end{table}

\begin{figure}[!ht]
	\begin{minipage}[b]{1.0\textwidth}
		\centering
		\includegraphics[height=5.5cm,width=9.6cm]{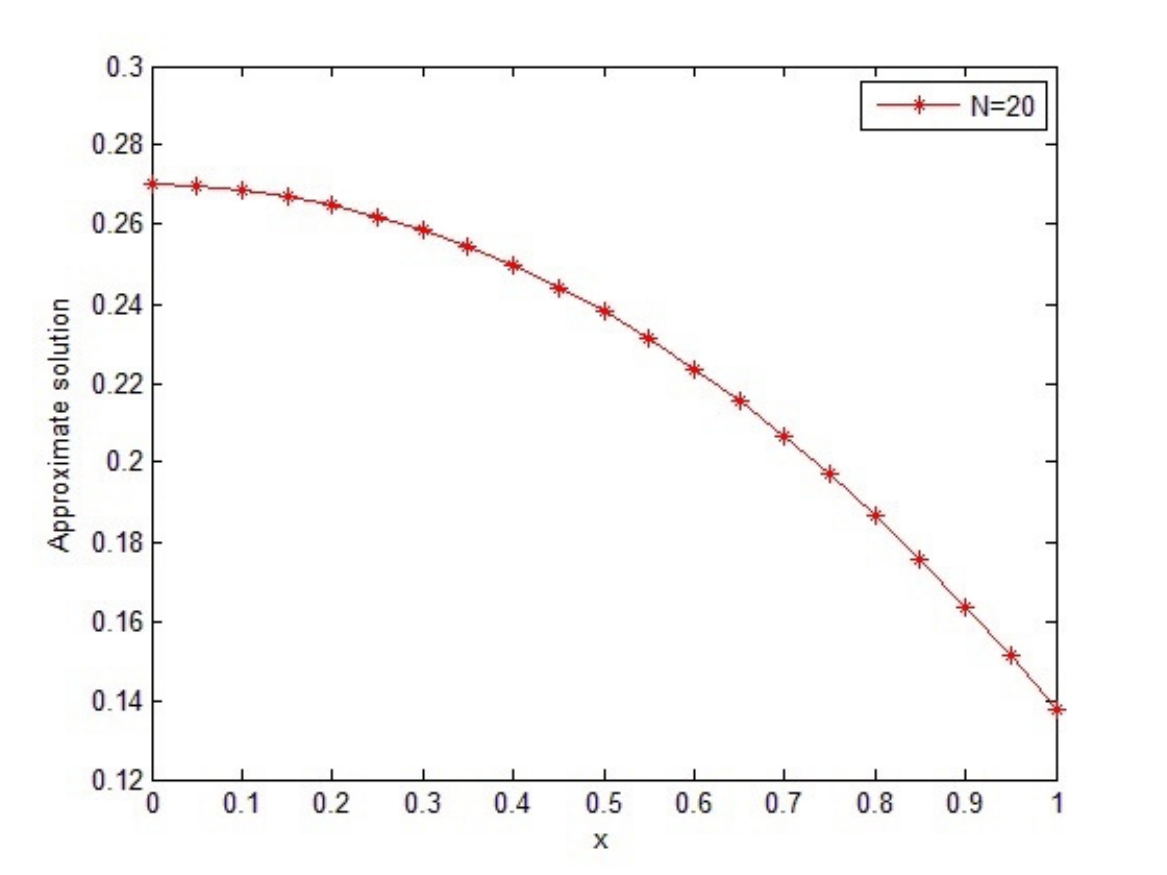}
		\caption{The graph of the approximate solution in Example 7 obtained by Method 2.}\label{figA}
	\end{minipage}
	\vfill
	\begin{minipage}[b]{1.0\textwidth}
		\centering
		\includegraphics[height=5.5cm,width=9cm]{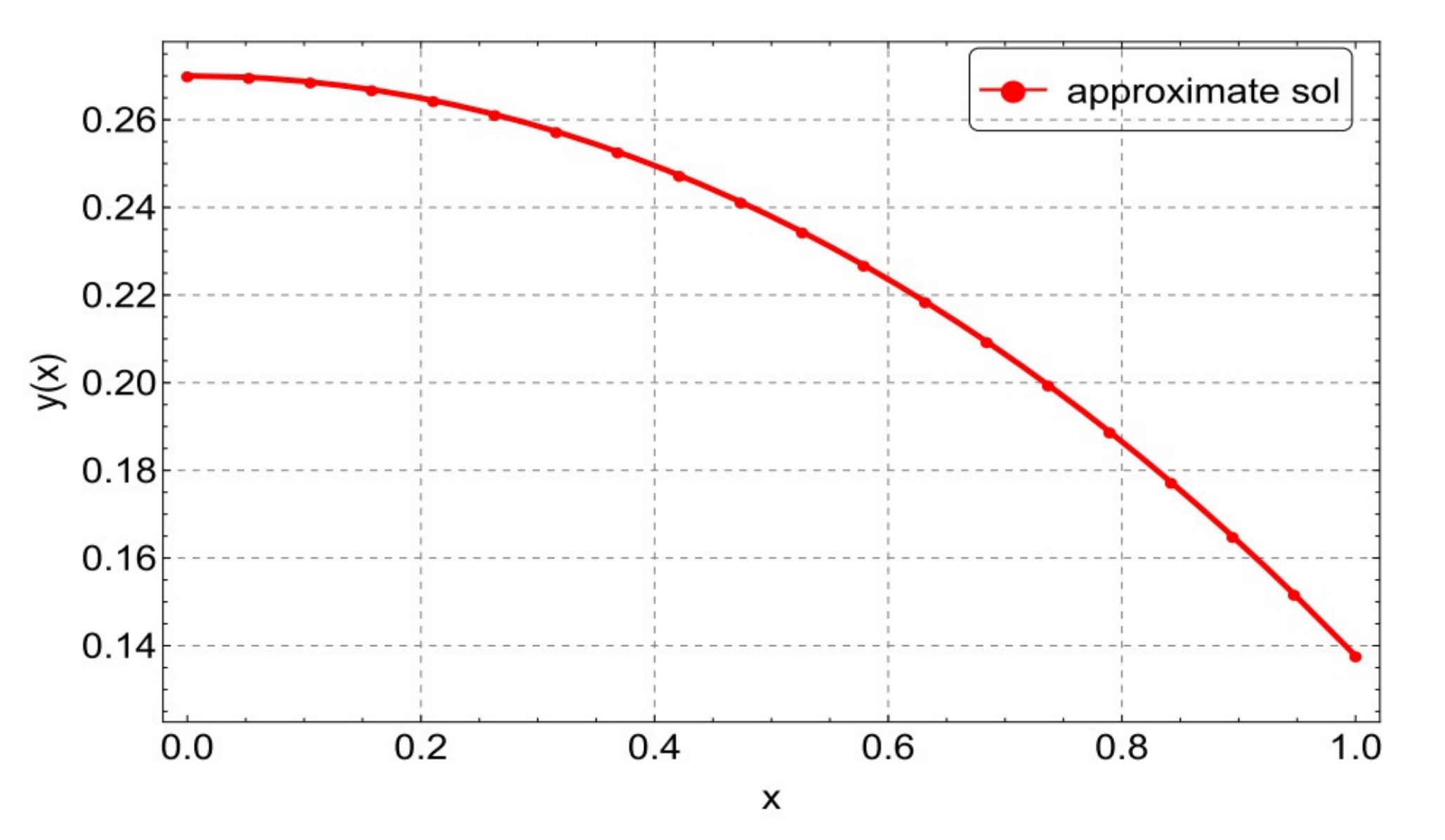}
		\caption{The graph of the approximate solution in Example 7 obtained in \cite{Malele}}\label{figM}
	\end{minipage}
\end{figure}

The graph of the approximate solution obtained for $N=20$ by our method and that of Malele in \cite{Malele} are depicted in Figures \ref{figA} and \ref{figM}, respectively.

\section{Conclusions and Future Perspectives}

This study has introduced a class of eighth-order numerical schemes for solving boundary value problems governed by the Lane-Emden equation. The proposed methods are derived from a continuous-level iterative formulation, wherein the integral representation of the solution is discretized using an enhanced trapezoidal rule augmented with correction terms from the Euler-Maclaurin formula. The schemes are straightforward to implement. A variety of numerical experiments confirm that the methods are highly accurate and superior when compared with recent state-of-the-art techniques.

Future work will aim to extend the current framework to address higher-order nonlinear singular models, such as third-order Emden-Fowler equations, which frequently arise in physical and chemical applications. While fourth-order methods for such equations have recently been proposed in the literature \cite{Alam1,Sahoo}, the development of higher-order schemes for these problems remains an open and promising direction.\\
{\bf Declarations
	Conflict of interest: } The authors declare that there is no Conflict of interest.


\begin{thebibliography}{99}

\bibitem {Chawla87} Chawla, M.M., Shivakumar, P.N.: On the existence of solutions of a class of singular nonlinear two-point boundary value problems. J. Comput. Appl. Math. {\bf 19} 379--388 (1987). https://doi.org/10.1016/0377-0427(87)90206-8

\bibitem {Pandey97} Pandey, R.K.: On a class of regular singular two point boundary value problems. J. Math. Anal. Appl. {\bf 208} 388--403 (1997). https://doi.org/10.1006/jmaa.1997.5320

\bibitem {Hoss} Hosseini, S.G., Abbasbandy, S.: Solution of Lane-Emden type equations by combination of the spectral method and Adomian decomposition method. Math. Probl. Eng. {\bf 2015} 1--10 (2015). https://doi.org/10.1155/2015/534754

\bibitem  {Liao} Liao, S.: A new analytic algorithm of Lane-Emden type equations. Appl. Math. Comput. {\bf 142}(1) 1--16 (2003). https://doi.org/10.1016/S0096-3003(02)00943-8

\bibitem {Singh}  Singh, R., Kumar, J.: An efficient numerical technique for the solution of nonlinear singular boundary value problems. Comput. Phys. Commun. {\bf 185}(4) 1282--1289 (2014). 
https://doi.org/10.1016/j.cpc.2014.01.002

\bibitem {Wazwaz}  Wazwaz, A.M.: A new algorithm for solving differential equations of Lane-Emden type. Appl. Math. Comput. {\bf 118}(2-3) 287--310 (2001). https://doi.org/10.1016/S0096-3003(99)00223-4

\bibitem {Ghor} Ghorbani, A., Bakherad, M.: A variational iteration method for solving nonlinear Lane-Emden problems. New Astron. {\bf 54} 1--6  (2017). https://doi.org/10.1016/j.newast.2016.12.004

\bibitem {Kanth} Ravi Kanth, A.S.V., Aruna, K.: He’s variational iteration method for treating nonlinear singular boundary value problems. Comput. Math. Appl. {\bf 60}(3) 821--829 (2010). https://doi.org/10.1016/j.camwa.2010.05.029

\bibitem {Wazwaz1} Wazwaz, A.M.: Variational iteration method for solving nonlinear singular boundary value problems arising in various physical models. Commun. Nonlinear Sci. Numer. Simul.
{\bf 16}(10) 3881--3886 (2011). https://doi.org/10.1016/j.cnsns.2011.02.026 

\bibitem {Khan} Khan, Y., Svoboda, Z., Šmarda, Z.: Solving certain classes of Lane-Emden type equations using the differential transformation method. Adv. Differ. Equ. {\bf 2012}(1) 1--11 (2012). https://doi.org/10.1186/1687-1847-2012-174

\bibitem {Singh3} Singh, O.P., Pandey, R.K., Singh, V.K.: An analytic algorithm of Lane-Emden type equations arising in astrophysics using modified homotopy analysis method. Comput. Phys. Commun. {\bf 180}(7) 1116--1124 (2009). https://doi.org/10.1016/j.cpc.2009.01.012

\bibitem {Chawla} Chawla, M.M., Mckee, S., Shaw, G.: Order $h^2$ method for singular two-point boundary value problem. BIT Numer. Math. {\bf 26} 318--326 (1986). https://doi.org/10.1007/BF01933711

\bibitem {Chawla1} Chawla, M.M., Subramanian, R., Sathi, H.L.: A fourth order method for a singular two-point boundary value problem. BIT Numer. Math. {\bf 28}(1) 88--97 (1988). https://doi.org/10.1007/BF01934697

\bibitem {Kumar} Kumar, M.: A difference scheme based on non-uniform mesh for singular two-point boundary value problems. Appl. Math. Comput. {\bf 136} 281--288 (2003). https://doi.org/10.1016/S0096-3003(02)00038-3

\bibitem {Verma} Verma, A.K., Kayenat, S.: On the convergence of Mickens’ type nonstandard finite difference schemes on Lane-Emden type equations. J. Math. Chem. {\bf 56} 1667--1706 (2018). https://doi.org/10.1007/s10910-018-0880-y

\bibitem {Malele} Malele, J., Dlamini, P., Simelane, S.: Solving Lane-Emden equations with boundary conditions of various types using high-order compact finite differences. Appl. Math. Sci. Eng. {\bf 31}(1) 2214303 (2023). https://doi.org/10.1080/27690911.2023.2214303

\bibitem {Roul3} Roul, P., Goura, V.P., Agarwal, R.: A compact finite difference method for a general class of nonlinear singular boundary value problems with Neumann and Robin boundary conditions. Appl. Math. Comput. {\bf 350} 283--304 (2019). https://doi.org/10.1016/j.amc.2019.01.001

\bibitem {Alam} Alam, M.P., Begum, T., Khan, A.: A high-order numerical algorithm for solving Lane-Emden equations with various types of boundary conditions. Comput. Appl. Math. {\bf{40}} 204 (2021). https://doi.org/10.1007/s40314-021-01591-7

\bibitem {Caglar} \c Ca\v glar, H., \c Ca\v glar, N., {\" O}zer, M.: B-spline solution of non-linear singular boundary value problems arising in physiology. Chaos Solitons Fract. {\bf 39}(3) 1232--1237  (2009). https://doi.org/10.1016/j.chaos.2007.06.007

\bibitem {Kada} Kadalbajoo, M.K., Kumar, V.: B-spline method for a class of singular two-point boundary value problems using optimal grid. Appl. Math. Comput. {\bf 188}(2) 1856--1869 (2007). https://doi.org/10.1016/j.amc.2006.11.050

\bibitem {Kanth1} Ravi Kanth, A.S.V.: Cubic spline polynomial for non-linear singular two-point boundary value problems. Appl. Math. Comput. {\bf 189} 2017--2022 (2007). https://doi.org/10.1016/j.amc.2007.01.002

\bibitem {Roul1} Roul, P., Thula, K., Agarwal, R.: Non-optimal fourth-order and optimal sixth-order B-spline collocation methods for Lane-Emden boundary value problems. Appl. Numer. Math. {\bf 145} 342--360 (2019). https://doi.org/10.1016/j.apnum.2019.05.004

\bibitem {Roul2} Roul, P., Thula, K., Goura, VMK P.: An optimal sixth-order quartic B-spline collocation method for solving Bratu-type and Lane-Emden–type problems. Math. Meth. Appl. Sci. {\bf 42}(8) 2613--2630 (2019). https://doi.org/10.1002/mma.5537

\bibitem {Roul4} Roul, P.: A high-order B-spline collocation method for solving a class of nonlinear singular boundary value problems. J. Math. Chem. {\bf 62} 1308--1322 (2024). https://doi.org/10.1007/s10910-024-01590-z

\bibitem {Dang2024-1} Dang, Q.A, Dang, Q.L., Ngo, T.K.Q.: Numerical methods of fourth, sixth and eighth orders convergence for solving third order nonlinear ODEs. Math. Comput. Simul. {\bf 221} 397--414 (2024). https://doi.org/10.1016/j.matcom.2024.03.018

\bibitem {Dang2024-2} Dang, Q.A, Nguyen, T.H., Vu, V.Q.: Construction of high order numerical methods for solving fourth order nonlinear boundary value problems. Numer. Algor. {\bf 99} 323--354 (2025). https://doi.org/10.1007/s11075-024-01879-9

\bibitem {Dang2025-1} Dang, Q.A, Nguyen, T.H., Vu, V.Q.: Eighth order numerical method for solving second order nonlinear BVPs and applications. J. Appl. Math. Comput.  {\bf 71} 3577-3600 (2025). https://doi.org/10.1007/s12190-025-02368-5

\bibitem {A-Long2021} Dang, Q.A, Dang, Q.L.: Simple numerical methods of second and third-order convergence for solving a fully third-order nonlinear boundary value problem. Numer. Algor. {\bf 87} 1479--1499 (2021). https://doi.org/10.1007/s11075-020-01016-2 

\bibitem {Singh1}  Singh, R., Nelakanti, G., Kumar, J.: Approximate solution of two-point boundary value problems using adomian decomposition method with Green function. Proc. Natl. Acad. Sci., India, Sect. A Phys. Sci. {\bf 85} 51--61 (2015). https://doi.org/10.1007/s40010-014-0170-4

\bibitem {Hild} Hildebrand, F.B.: Introduction to numerical analysis. 2nd edition, McGraw Hill, New York (1974)

\bibitem {Dickey} Dickey, R.: Rotationally symmetric solutions for shallow membrane caps. Quart. Appl. Math. {\bf 47}(3) 571--581 (1989). https://doi.org/10.1090/qam/1012280

\bibitem {Singh2} Sahoo, N., Singh, R., Kanaujiya, A., Cattani, C.: A new fourth-order compact finite difference method for solving Lane-Emden-Fowler type singular boundary value problems. J. Comput. Sci. {\bf 83} 102474 (2024). https://doi.org/10.1016/j.jocs.2024.102474



\bibitem {Lions} Lyons, M.E.: Charge percolation in electroactive polymers.  1--235. In: Electroactive Polymer Electrochemistry: Part 1: Fundamentals, Plenum Press, New York (1994)

\bibitem  {Duggan} Duggan, R., Goodman, A.: Pointwise bounds for a nonlinear heat conduction model of the human head. Bull. Math. Biol. {\bf 48}(2) 229--236 (1986). https://doi.org/10.1016/S0092-8240(86)80009-X




\bibitem {Alam1} Alam, M.P., Khan, A.: An efficient collocation algorithm for third order non-linear Emden-Fowler equation. Soft Comput. {\bf 29} 1767--1788 (2025). https://doi.org/10.1007/s00500-025-10431-3

\bibitem {Sahoo} Sahoo, N., Singh, R.: Compact finite difference schemes and error estimation for third-order Emden-Fowler equations. Numer. Algor. (2025) https://doi.org/10.1007/s11075-025-02025-9

\end{thebibliography}
\end{document}